%% file: main_blank.tex
\documentclass{article}

\input{header.tex}

\usepackage[a4paper, margin=1in]{geometry}
\usepackage{biblatex}
\usepackage{orcidlink}
\author{Alexander Sikorski\,\orcidlink{0000-0001-9051-650X}\thanks{Corresponding author, \tt{sikorski@zib.de}}, 
Enric Ribera Borrell\,\orcidlink{0000-0002-3897-5984},
Marcus Weber\,\orcidlink{0000-0003-3939-410X} \\
Zuse Institute Berlin}

\begin{document}
\title{Learning Koopman eigenfunctions of stochastic diffusions with optimal importance
sampling and ISOKANN\thanks{The following article has been submitted to the Journal of Mathematical Physics. After it is published, it will be found at \href{https://aip.scitation.org/journal/jmp}{Link}.}} 
\maketitle
\begin{abstract} \theabstract \end{abstract}

\input{content.tex}


\end{document}

%% file: header.tex
\usepackage[utf8]{inputenc}

\usepackage{mathtools}
\usepackage{amsmath}
\usepackage{amsthm}
\usepackage{amsfonts}
\usepackage{amssymb}

\usepackage{algorithm}
\usepackage[indLines=true]{algpseudocodex}
\algrenewcommand\algorithmicrequire{\textbf{Input: \,\,}}
\algrenewcommand\algorithmicensure{\textbf{Output:}}
\usepackage{caption}
\usepackage{subcaption}

\usepackage{hyperref}
\hypersetup{
     colorlinks=true,
     linkcolor=blue,
     anchorcolor=black,
     citecolor=blue,
     filecolor=black,      
     urlcolor=black,
}

\usepackage{cleveref}
\crefname{equation}{Eq.}{Eqs.}
\crefname{section}{Sec.}{Secs.}
\Crefname{section}{Section}{Sections}

\usepackage{xspace}

\newcommand{\RR}{\mathbb{R}} 
\newcommand{\EE}{\mathbf{E}} 
\newcommand{\dd}{\mathrm{d}} 
\newcommand{\KK}{\mathbf{K}} 

\newcommand{\one}{\mathbf{1}} 
\newcommand{\ef}{v} 

\newcommand{\pccap}{\texorpdfstring{\texttt{PCCA$^+$}}{PCCA+}\xspace}
\newcommand{\XX}{\mathbf{X} } 
\newcommand{\isokann}{\texorpdfstring{\texttt{ISOKANN}}{ISOKANN}\xspace}
\newcommand{\isokannI}{\texorpdfstring{\texttt{1D-ISOKANN}}{1D-ISOKANN}\xspace}
\newcommand{\Id}{\text{Id}\xspace}

\newtheorem{theorem}{Theorem}
\newtheorem{proposition}{Proposition}
\newtheorem{corollary}{Corollary}[theorem]

\newcommand{\theacknowledge}{
We thank Luca Donati and Luzie Helfmann for their support in proofreading and insightful discussions.
This research has been funded by Deutsche Forschungsgemeinschaft (DFG) through grant CRC 1114 "Scaling Cascades in Complex Systems", Project Number 235221301, Project A05 "Probing scales in equilibrated systems by optimal nonequilibrium forcing".
}

\newcommand{\theabstract}{
For stochastic diffusion processes the dominant eigenfunctions of the corresponding Koopman operator contain important information about the slow-scale dynamics, that is, about the location and frequency of rare events.
In this article, we reformulate the eigenproblem in terms of $\chi$-functions in the \isokann framework and discuss how optimal control and importance sampling allows for zero variance sampling of these functions. We provide a new formulation of the \isokann algorithm allowing for a proof of convergence and incorporate the optimal control result to obtain an adaptive iterative algorithm alternating between importance sampling and $\chi$-function approximation.
We demonstrate the usage of our proposed method in experiments increasing the approximation accuracy by several orders of magnitude.

}

\newcommand{\thecontent}{
\input{content.tex}
}

%% file: content.tex
\newcommand{\Sbar}{\bar S}
\newcommand{\koopvar}{\kappa}
\newcommand{\Niter}{N}
\newcommand{\Nx}{M}
\newcommand{\Nmc}{K}
\newcommand{\Nlearn}{L}

\newcommand{\lag}{T}
\newcommand{\KKt}{\KK^\lag}
\newcommand{\chieq}{\emph{$\chi$-equation}}
\newcommand{\chifunctions}{$\chi$-functions\xspace}

\section{Introduction}

Many real-world stochastic processes contain rare events, for example folding and binding events in molecular systems. The analysis of the frequency and mechanism of these events often takes operators associated with the process and their dominant invariant subspaces into account  \cite{Bovier2004, Huisinga2021, Nielsen2016, Rabben2020isokann, Zhang2022Koopmanfr}.
Usually this type of analysis leads to a kind of chicken and egg problem. In order to compute the dominant invariant subspace of the Koopman operator of the process, one has to somehow also ``sample those events which are usually rare''. However, in order to know how to generate a bias on the process for observing these events, one would need information from the dominant invariant subspace of the Koopman operator.

The key idea of this article is to take an iterative algorithm which approximates the dominant eigenfunctions of the operator and to use the intermediate approximations for sampling along the reaction paths and generating an optimal bias to observe the relevant events.  The algorithm for approximating eigenfunctions is \isokann\footnote{An acronym for ``Invariant subspaces of Koopman operators with artificial neural networks''.} \cite{Rabben2020isokann} and the bias is computed according to the theory of optimal importance sampling and the change of path measures \cite{Hartmann2018optimalcontrol}.

Briefly speaking, \isokann can be thought of as an Arnoldi-like method using neural networks as function representations and replacing the subspace projections by a transformation suitable to its application on neural networks. It does not compute the eigenfunctions themselves but rather the so called \chifunctions, which span an invariant subspace of the Koopman operator and can be used to reconstruct the eigenfunctions. Furthermore the \chifunctions themselves allow for an interpretation as reaction coordinates indicating the locations of rare events and their reaction paths \cite{Ernst2019, SchuetteKlusHartmann2022}.

 Using this interpretation, the \chifunctions obtained during previous iterations can be used to adapt the sample locations for further iterations, e.g. by \emph{$\chi$-stratified sampling} ({\cref{sec:computation}), and thus providing better global coverage and facilitating exploration.
The theory of optimal importance sampling on the other hand allows to exploit the information locally by decreasing the variance of the samples required to approximate the action of the Koopman operator.
This in turn allows \isokann to arrive at results either more quickly or more precisely.

It is due to this feedback loop, coupling local and global information, together with its representation of the \chifunctions as neural networks that we expect \isokann to perform well even for complex high-dimensional systems.
While these general ideas are far from being formalized yet, in this article we will try to construct the basic building blocks in a way amenable for their future analysis.


In the following, we start by recalling basic knowledge about the eigenfunctions of Koopman operators (\cref{sec:koopman}) and provide an abstract formulation of the \isokann problem in terms of \chifunctions  (\cref{sec:isokann}) encoding the invariant subspaces of the Koopman operator.
This abstract formulation is incomplete without the choice of an adequate transformation which we discuss in \Cref{sec:transformationS} and follow it up with an explicit choice leading to \isokannI (\cref{sec:isokannI}) which corresponds to the classical \isokann algorithm \cite{Rabben2020isokann}. After proving convergence of \isokannI we show how to reconstruct the eigenfunctions of the Koopman operator from the \chifunctions in \Cref{sec:eigenfunction}.
We then conclude \Cref{sec:isotheory} by providing an algorithmic description in form of \Cref{alg:isokann} and discussing the actual sampling and computation procedure (\cref{sec:computation}).
\Cref{sec:optcontrol} starts with an introduction to the theory of optimal importance sampling of classical random variables (\cref{sec:oisRV}) before reciting the result for the optimal sampling of path observables for diffusion processes with the Girsanov re-weighting in \Cref{sec:oisdiff}. After showing how to apply this result to obtain a zero variance sampler for the evaluation of the Koopman operator (\cref{sec:oiskoop}) we discuss how to integrate the result into the \isokann framework (\cref{sec:oisiso}).
In the end (\cref{sec:oisexample}) we apply the \isokannI algorithm to a one-dimensional double-well potential and observe the improvement of the accuracy of the controlled versus the uncontrolled case.
 
\section{\isokann Theory}\label{sec:isotheory}

Before introducing our key idea in the next chapter, we will in this chapter recall the basics of Koopman operator theory and summarize the \isokann method while also complementing it by some new results. In particular we provide a new dimension-agnostic formulation \eqref{eqn:isokann} which naturally leads to a possible extension of \isokann to higher dimensions and show how the classical \isokannI can be seen as a special case of this formulation and prove convergence of the algorithm to a $\chi$ function (\Cref{thm:convergence}). After showing how to reconstruct the dominant eigenfunctions from the \isokann result (\Cref{prop:eigenfunction}) we conclude this section with a discussion of the actual implementation of \isokann, suggesting a new adaptive sampling scheme.

\subsection{The Koopman operator}\label{sec:koopman}

\begin{figure}
\begin{minipage}[t]{.485\textwidth}
  \centering
  \includegraphics[width=\linewidth]{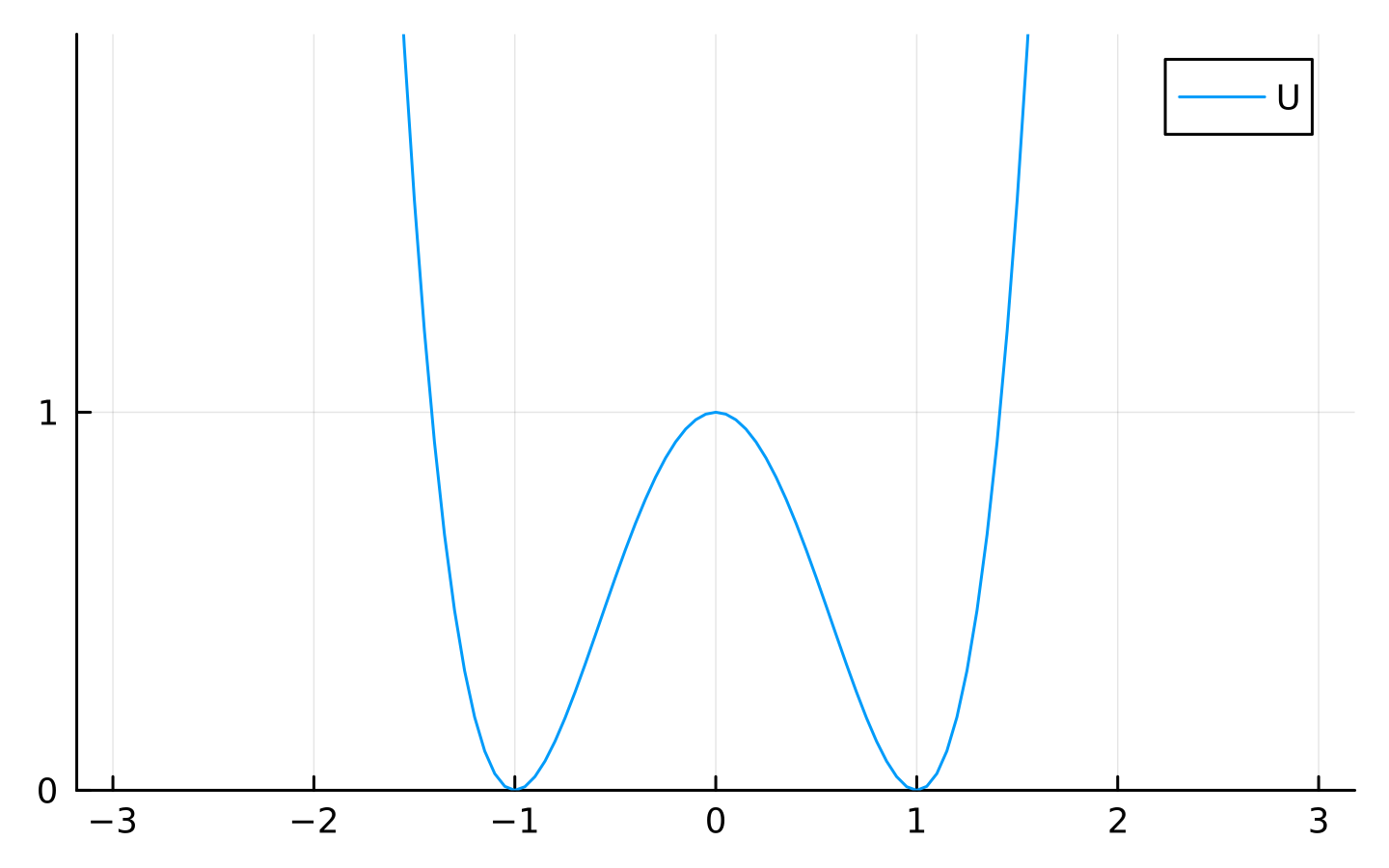}
  \caption{Potential function $U$ of the double well with two metastable regions separated by a potential barrier.}
  \label{fig:potential}
\end{minipage}%
\hfill
\begin{minipage}[t]{.485\textwidth}
  \centering
  \includegraphics[width=\linewidth]{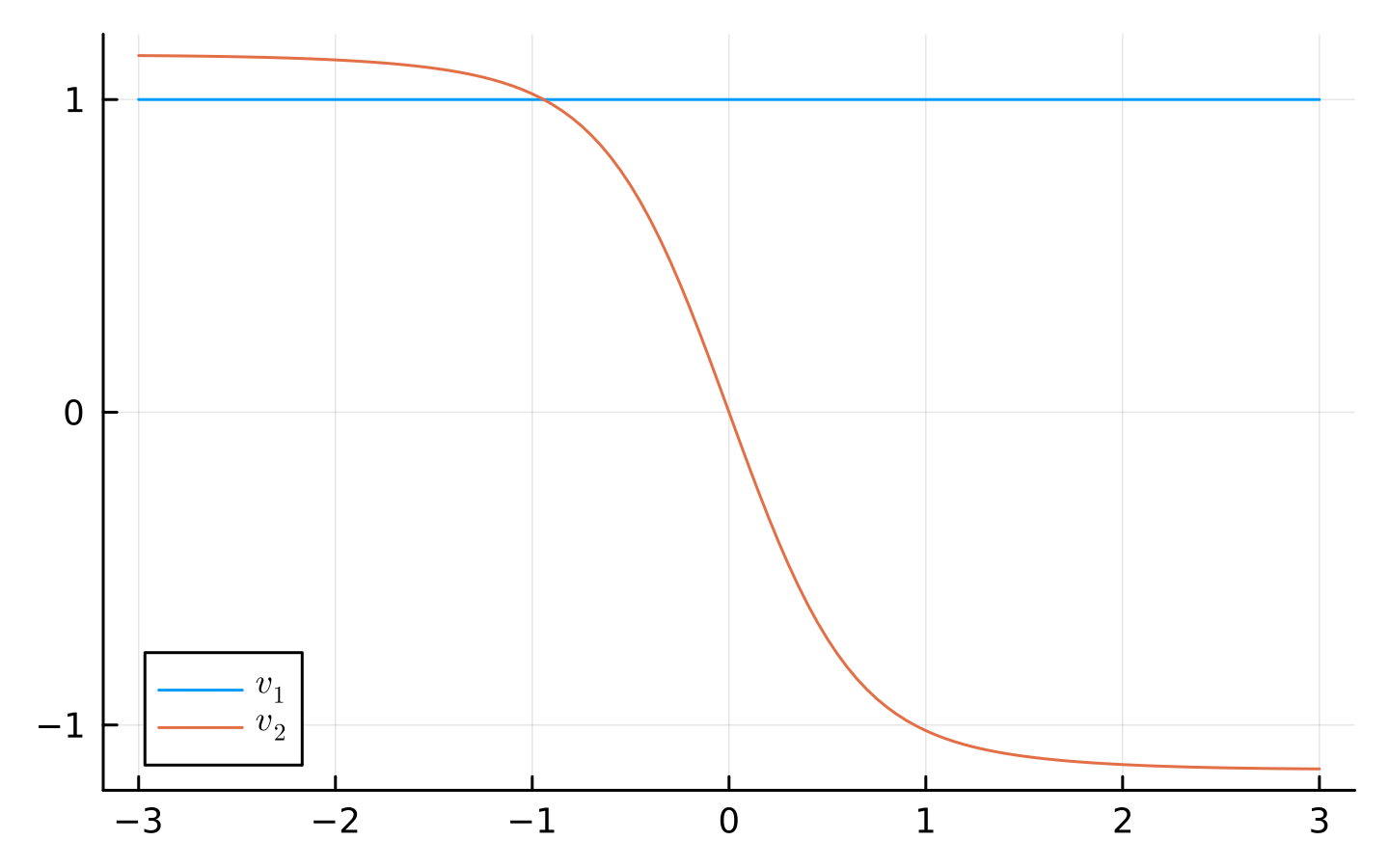}
  \caption{The two dominant eigenfunctions $v_1, v_2$ of the Koopman operator for the double well potential.}
  \label{fig:eigenfun}
\end{minipage}
\end{figure}

Although our approach can be generalized to non-reversible stochastic processes (by shifting the focus from eigenfunctions to invariant subspaces), for simplicity we will restrict our explanations to the reversible case. More precisely, we will investigate potential-driven diffusion processes  $X = (X_t)_{t \geq 0}$ of the form
\begin{equation}
\label{eq: not controlled sde}
\dd X_t = b(X_t) \dd t + \sigma \dd B_t,
\end{equation}
%
taking values in the state space $\XX =  \RR^n$  with constant diffusion term $\sigma \in \RR^{n \times n}$ and force-field $b: \XX \rightarrow \RR^n$ given by the gradient of a smooth potential $b = -\nabla U$, $U:\XX \rightarrow \RR$
\footnote{
We choose $b$ to be a gradient field for the process $X_t$ to be reversible and hence admit a real eigendecomposition. Our principal results do not require reversibility when arguing in terms of invariant subspaces instead of eigenfunctions. However, for the sake of simplicity we here consider reversible systems only.}.
$B$ is a $n$-dimensional Brownian motion. 

The Koopman operator for lag time $\lag$, 
$\KKt:L^\infty(\XX) \rightarrow L^\infty(\XX)$, 
applied to a function 
$f \in L^\infty(\XX)$
is defined by its pointwise evaluation via
\begin{align}
    \left(\KKt f \right) (x) = \EE \left [f(X_\lag) \mid X_0 = x\right]
\end{align}
i.e. the expectation value of $f$ at time $\lag$ when starting the system in $X_0 = x$. Recall that since the process is time-homogeneous the Koopman operator just depends on the lag time for any start time $t \geq 0$
\[
\EE \left [f(X_{t+\lag}) \mid X_t = x\right] = \left(\KKt f \right) (x) .
\]

The eigenfunctions $\ef_i \in L^\infty(\XX)$ of $\KKt$ satisfy for all lag times $\lag \geq 0$
\begin{align}
    \KKt \ef_i &= \lambda_{i}(\lag) \,\ef_i, && i=1,2,... \\
    \lambda_i(\lag) &= \exp(\lag q_i), && 0 = q_1>q_2\ge...
\end{align}
with time-dependent eigenvalues $\lambda_i(\lag)$, exponential in time with rates $q_i$ (which in turn are the eigenvalues of the corresponding infinitesimal generator of $\KKt$).
In the following we will refer to $\ef_1, \dots, \ef_d$ as the $d$ \emph{dominant eigenfunctions} and call $\ef_1\equiv 1$ the \emph{trivial eigenfunction}. When clear from the context or of no importance we will omit the lag time $\lag$ and simply speak of the Koopman operator $\KK$.

The dominant eigenfunctions are of particular interest as they decay the slowest and hence dominate the long time behavior of the system. The number $d$ of eigenfunctions of interest depends on the time scales of the system and is usually chosen up to a spectral gap.

There exist different approaches to estimate the eigenfunctions. 
Many depend on the discretization of the state-space into cells leading to a matrix representation of $\KK$.
A classical method is starting trajectories in each such cell and counting how many of them end up in a certain cell. This sample-driven method can be interpreted as an approximate Ulam/Galerkin discretization onto indicator functions of the cells \cite{Klus2016}. The eigenfunctions of $\KK$ are then approximated by the eigenfunctions of its (dense) matrix approximation.
Unfortunately this scheme breaks down in high dimensions as the number of cells in a structured grid increases exponentially.

One possible remedy to this problem is posed by the Square-Root-Approximation method (\texttt{SQRA})\cite{Donati2021markov}. It approximates the infinitesimal generator of the Koopman operator by a finite volume approximation  where the volumes are implicitly defined by the Voronoi tesselation induced by some sample points. Using only evaluations of the potential $U$ at these samples it can be understood as a semi-parametric method resulting in a sparse matrix representation, which in turn can be used for the computation of the eigenfunctions. 

All these classical approaches however depend on a discretization {\em before} solving the eigenproblem. As an alternative we will now summarize \isokann, a recent matrix-free approach learning a linear combination of the eigenfunctions by neural networks.

\subsection{\isokann -- Computing the dominant eigenspace}\label{sec:isokann}

The \isokann algorithm \cite{Rabben2020isokann}
uses a nonparametric representation in the form of a neural network in order to learn the dominant invariant subspace by interleaving an Arnoldi-like power iteration with the approximation of the Koopman operator by Monte-Carlo simulations.

To this end it will be useful to reformulate the eigenproblem in terms of the so called \emph{$\chi$-function} $\chi = (\chi_i)_{i=1}^d: \XX \rightarrow \RR^d$  with
$0 \leq \chi_i\leq 1, \ \sum_i \chi_i=1$ satisfying the \chieq, 
\begin{equation}\label{eq:isofix}
    \chi = S \KK \chi,
\end{equation}
for an appropriately chosen matrix $S:\RR^d\rightarrow\RR^d$, which will be specified in the next section. The core idea here is that the components $\chi_i$ of $\chi$ span an invariant subspace of $\KK$, or more specifically they consists of linear combinations of eigenfunctions of $\KK$. Furthermore the action of $\KK$ on this space is then explicitly given by the matrix $S^{-1}$. 
This definition of $\chi$ functions does not only allow us to reconstruct the eigendecomposition of $\KK$ (\Cref{prop:eigenfunction}) but also allows for an interpretation as macrostates in the form of \emph{fuzzy memberships} to $d$ different sets \cite{Roblitz2013fuzzy} and leads to a direct characterization of rare transitions in form of their holding times, reaction rates, exit paths etc. \cite{Ernst2019,SchuetteKlusHartmann2022}.

Formally \isokann then approximates the $\chi$-function representing the dominant eigenspace by an iterative sequence of approximations 
$\chi_n: \XX \rightarrow \mathbb{R}^d$ satisfying the \emph{\isokann equation}, 
\begin{equation}\label{eqn:isokann}
    \chi_{n+1}(x) = S_n \KKt \chi_n(x) = S_n \EE \left [\chi_n(X_\lag) \mid X_0 = x\right],
\end{equation}
with $S_n=S_n(\KKt\chi_n)$ being a linear map depending on the previous iteration and constructed in such a way as to provide convergence to a desired form of $S$ in \eqref{eq:isofix}.


Similar to the Arnoldi-iteration, the iterative application of the Koopman operator leads to a decay of the eigenfunctions that is exponential in their eigenvalue.
In the following section we will discuss how to construct $S_n$ as to compensate this decay just so that the dominant components will prevail whilst the non-dominant ones vanish. In the limit, $\chi_n$ spans the dominant subspace (\Cref{thm:convergence})

\begin{align}\label{eqn:convergence}
    \text{span} \{(\chi_n)_1, ..., (\chi_n)_d\} \overset{n\rightarrow\infty}{\rightarrow} \text{span}\{\ef_1, ..., \ef_{d}\},
\end{align}
with $v_i$ denoting the dominant eigenfunctions.
This in turn allows us to learn the linear action of $\KK$ on that subspace leading to
\begin{align}
    \chi_n\rightarrow\chi \text{ and } S_n \rightarrow S.
\end{align}

Let us note that this presentation of \isokann differs from the original variant~\cite{Rabben2020isokann} in that the iteration in the form of \eqref{eqn:isokann} allows general $d$-dimensional valued $\chi$ functions.
However we will see that for the case of $d=2$ this is equivalent to the original variant, which by slight abuse of notation we will refer to as \isokannI.

\paragraph{Illustrative example}
To illustrate this better, let us consider the following example given by the SDE \eqref{eq: not controlled sde} with a double well potential 
\begin{equation}
\label{eqn:doublewell}
    U(x) = (x^2 - 1)^2.
\end{equation}
The double well potential with two minima is shown in \Cref{fig:potential}, the two dominant eigenfunctions are given in \Cref{fig:eigenfun}. In \Cref{fig:chifun}, we show the resulting two $\chi$-functions, that are a linear combination of the eigenfunctions, but allow us to grasp the slow time-scale dynamics better. In particular, the two $\chi$-functions can be interpreted as membership functions of the two wells of the potential. Their exit paths are charachterized by the gradients $-\nabla \chi$ and their holding probabilties or exit rates can be computed from the associated matrix $S$ in \eqref{eq:isofix} (see \cite{Ernst2019,SchuetteKlusHartmann2022}).

\subsection{Choice of the transformation \texorpdfstring{$S$}{S}}
\label{sec:transformationS}

\begin{figure}
\begin{minipage}[t]{.485\textwidth}
  \centering
  \includegraphics[width=\linewidth]{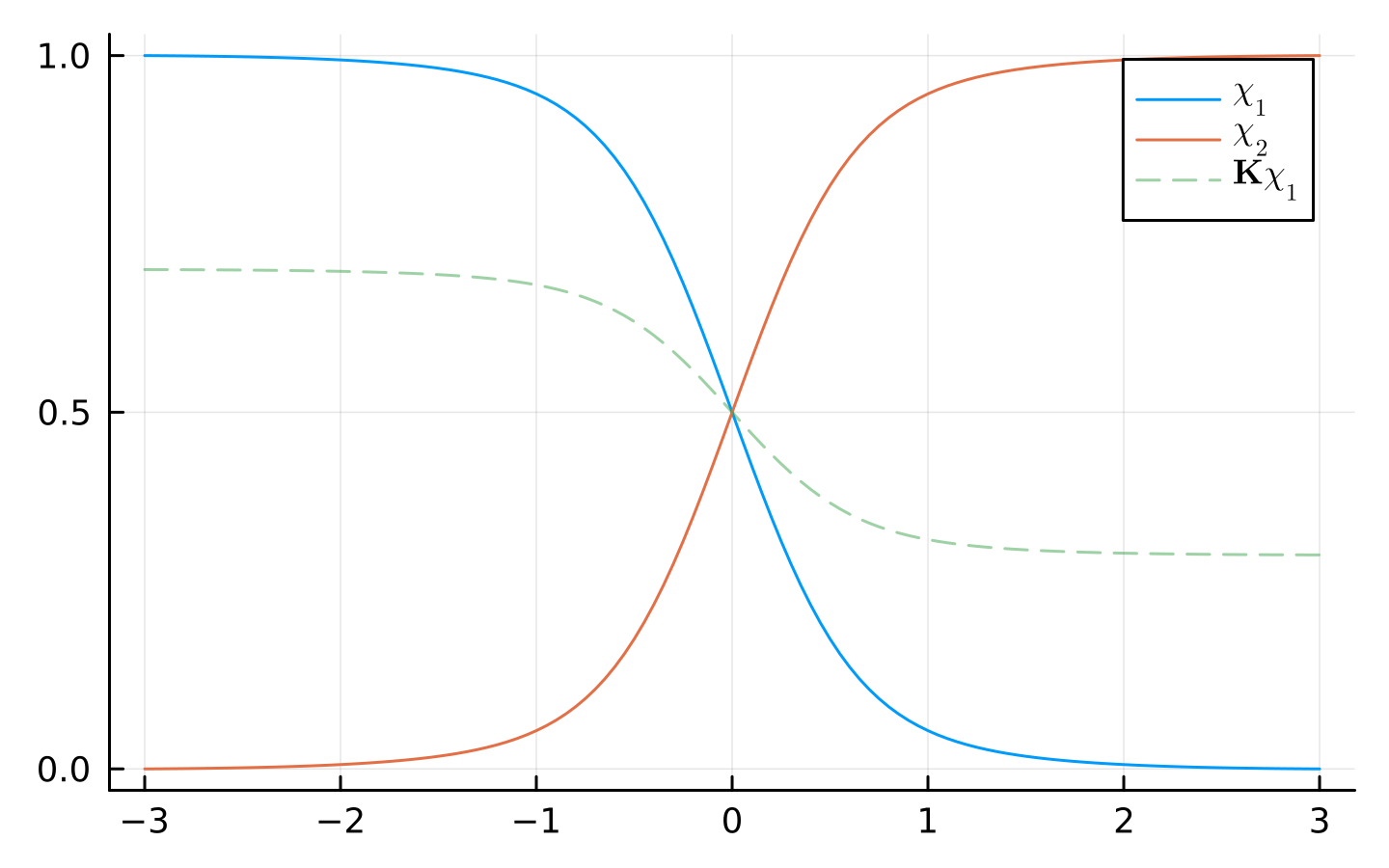}
  \caption{Both components of the two-dimensional $\chi$-function, which are linear combinations of the eigenfunctions (\Cref{fig:eigenfun}).
  In this case the application of $K$, which in general is linear in $\chi$, corresponds to a shift-scale (see \Cref{sec:isokannI})}
  \label{fig:chifun}
\end{minipage}%
\hfill
\begin{minipage}[t]{.485\textwidth}
  \centering
  \includegraphics[width=\linewidth]{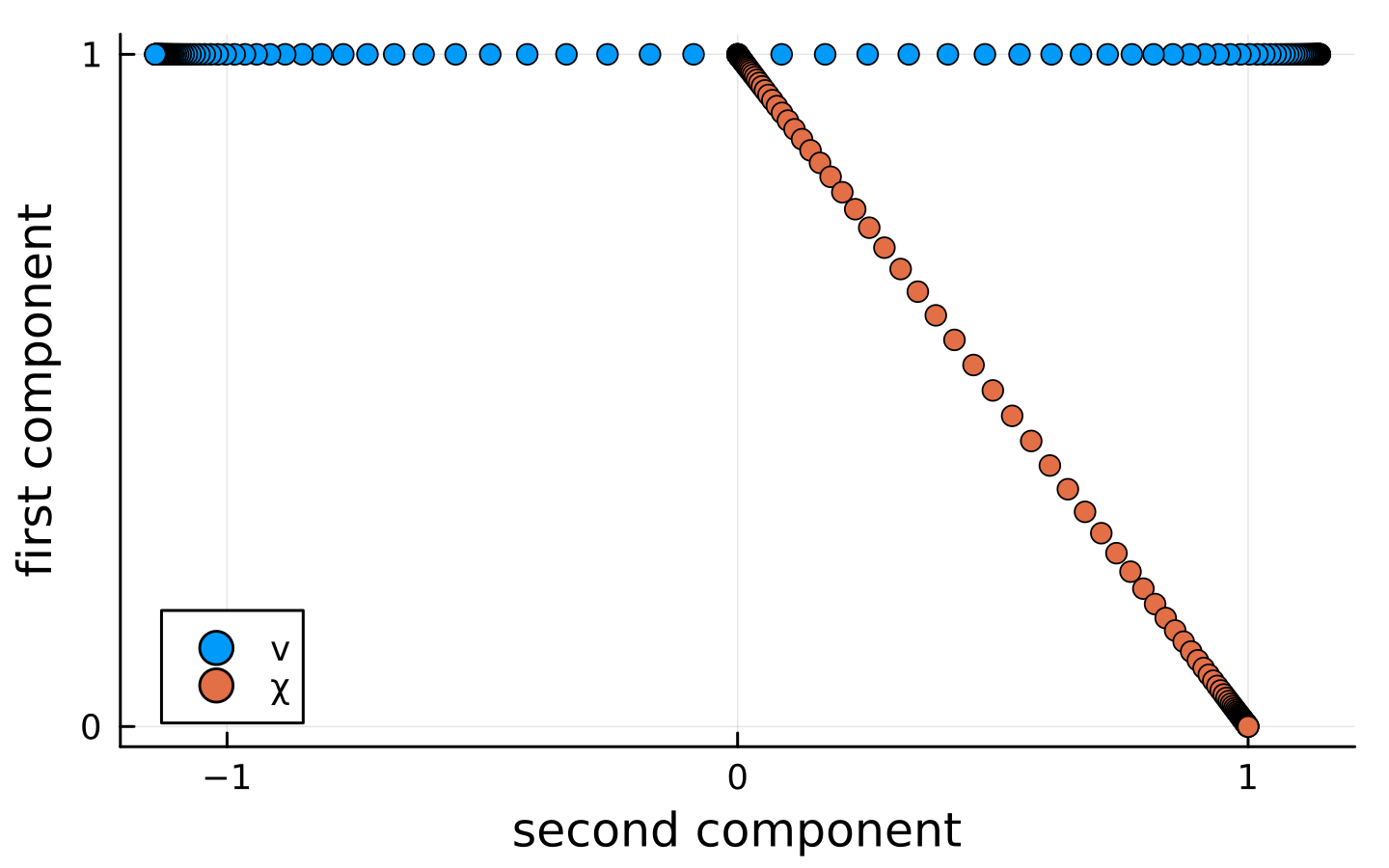}
  \caption{Scatter plot depicting $v_1$ against $v_2$ and $\chi_1$ against $\chi_2$ on a uniform grid over $\XX$. \pccap constructs the linear map $S$ mapping $v$ onto the unit-simplex~$\chi$. }
  \label{fig:compscat}
\end{minipage}
\end{figure}

In the previous description we were talking about an "appropriately chosen" linear transformation $S$ or a sequence of $S_n$. We will now explain the role that $S$ plays in \isokann and suggest a way of determining suitable $S$.

In principle the map $S$ can be chosen arbitrarily such that \eqref{eq:isofix} has a solution in $\chi$. However to obtain a useful result and for the \isokann iterations to converge with the corresponding choice of $S_n$ we require some specific properties.

\newcommand{\arnoldimeth}{\emph{Arnoldi method}}
In order to understand the role of $S$ let us for now think about the classical \arnoldimeth to find the dominant eigenfunctions of a matrix.
In principle the Arnoldi method also takes the form of \eqref{alg:isokann}, where the application of $S_n$ corresponds to a \emph{Gram-Schmidt orthonormalization}.
Here the orthogonalization ensures that the leading eigenvectors are projected out from the subsequent ones and the following normalization step then ensures that the eigenvectors do not decay over multiple iterations.

Note here that while the Gram-Schmidt orthonormalization on its own is a nonlinear procedure, the resulting action on a given set of input vectors can be expressed as a linear map, i.e. if $\omega$ is the orthonormalization procedure, for each input matrix $k$ we can find a linear map $\Omega(k)$ (depending on $k$) such that
\begin{equation}
    \label{eqn:linearnonlinear}
    \omega(k) = \Omega(k) \, k.
\end{equation}
It is in exactly this way that we understand $S$ as a linear map computed non-linearly on the data in \cref{eqn:isokann}.

For \isokann we want $S$ to fulfill the same role: The goal of each linear transformation $S_n$ is to counteract the decay of the dominant eigenfunction components contained in $\chi_n$ after application of the Koopman operator $\KK$.

However, we are interested in linear combinations of the eigenfunctions over a \emph{continuous space} represented by a neural network. In this setting orthonormalization is a hard problem, involving integration over the whole state space
\footnote{One might resort to restricting the $\chi$ functions to a finite number of fixed points and orthonormalize wrt. the resulting vectors. We expect this to work fine as long as including those points where the eigenfunction differ strongly, i.e. the individual metastabilites. However, since these were not allowed to change and are usually not known a priori they do not lend themselves to an adaptive scheme like \isokann does.}
.

On the other hand, orthonormality is a strong assumption which is not necessarily required. If we manage to choose $S$ such that it amplifies the first $d$ eigenfunctions such that they stay bounded away from zero, we will obtain a representation of the dominant subspace, which in sequence allows for the reconstruction of eigenfunctions (c.f.  \Cref{sec:eigenfunction}).

We now motivate a heuristic approach, based on the \pccap algorithm \cite{Roblitz2013fuzzy}, to construct a transformation $S$ and will prove its convergence in the 1D-setting.

Let us shortly summarize the idea of the \pccap methodology.
In the context of metastable systems the state-space regions where the individual eigenfunctions become extremal are representative for the respective metastabilities of that system. In practice, plotting the state space $\XX$ over the respective eigenfunction components $\ef_i$ the resulting set often resembles a simplex.
\pccap can be understood as a method to identify this simplex structure and construct the linear transformation (in the space of eigenfunction components) mapping this set into the unit-simplex (see also \Cref{fig:compscat}), such that the image becomes "as big as possible" (specified by an optimization problem).

We can similarly imagine this picture for the $\chi$ functions:
After sufficient time propagation (or power iterations), the $\chi$ functions are mainly composed of the dominant eigenfunctions. Plotting $\XX$ over the $\chi$ components thus will be close to the $\XX$ over $\ef$ plot above modulo a linear transformation/a change of basis.
We can thus apply \pccap onto our intermediate $\KK \chi_n$ in order to find a transformation such that the resulting $\chi_{n+1}$ fills this unit-simplex. This inhibits the exponential decay of the non-trivial eigenfunctions by maintaining a set of $d$ linearly independent components. Since the dominant eigenfunctions decay slower then the following non-dominant ones, they will dominate the behavior and prevail for $n\rightarrow\infty$.

Even though preliminary results have shown that using \pccap (with minor modifications\footnote{In order for this to work \pccap has to be stable with respect to permutations, i.e. one has to ensure that the $S_n$ do not suddenly change the orientation, which would prevent convergence.}) to construct $S$ in higher dimensions works fine, the focus of this paper is on the optimal control so we will reserve a more detailed report to future work and hence give a proof only for the simpler case of a single time-scale in the next section.

\subsection{\isokannI with explicit \pccap}
\label{sec:isokannI}
In the case where one is interested merely in the first non-trivial eigenfunction, $\ef_2$, the \pccap solution can be computed explicitly, which in turn allows us to prove convergence of \isokann for $d=2$ and explicitly specify $S$.

Note that because of $\chi_1 + \chi_2 =1$, the desired solution $\chi:\XX\rightarrow\RR^2$ is fully determined by its first component $\bar\chi := \chi_1:\XX\rightarrow\RR$ alone,
\begin{equation} \label{eqn:chibar}
\chi = \begin{pmatrix}\bar\chi \\1-\bar\chi\end{pmatrix}.
\end{equation}

This representation allows us to solve the two dimensional \isokann problem with just one scalar function, approximated by a series of scalar neural networks $\bar\chi_n:\XX\rightarrow\RR$. This scalar representation is the approach taken in \cite{Rabben2020isokann} and in the following we will refer to it as \isokannI\footnote{We see this abuse of notation justified as \isokann for $d=1$ would otherwise only denote the trivial solution $\chi\equiv1$}.

Let us start by constructing the explicit \pccap solution.
Recall that $\ef_1 \equiv 1$ and note that $\ef_2(\XX) = \left[ \min(\ef_2), \max(\ef_2) \right]$. Thus the image of $\XX$ forms a line segment, i.e. a 1-dimensional simplex, in the $\ef_1$-$\ef_2$-plane.
\pccap then constructs the unique map $S$ which maps this simplex onto the unit simplex $\left\{(x,y) \mid x+y=1, x>0, y>0\right\}$ (see \Cref{fig:compscat}).

To this end, let us introduce the map $\Sbar$ for bounded continuous functions $\koopvar\in C(\XX)$
\newcommand{\Sarg}{\koopvar}
\begin{equation}\label{eqn:shiftscale}
    \Sbar(\Sarg) = \frac
    {\Sarg - \min(\Sarg)}
    {\max(\Sarg) -\min(\Sarg)}
\end{equation}
such that $\Sbar(\koopvar): \XX \twoheadrightarrow [0,1]$ maps surjectively onto the unit interval.
Even though $\Sbar$ (consisting of a shift and a scale) is only affine-linear in its argument it can be seen as a linear map on the 1D subspaces $\{(1,\koopvar) \mid \koopvar\in C(\XX)\})$ or similarly $\{(\koopvar,1-\koopvar) \mid \koopvar\in C(\XX)\}$ and indeed is the action of the \pccap solution on a single component, i.e. there exists a matrix $S \in \RR^{2\times 2}$ such that it satisfies
\begin{equation} \label{eqn:fullS}
S\cdot\begin{pmatrix}
     \koopvar\\  
     1-\koopvar 
 \end{pmatrix} = \begin{pmatrix}
     \Sbar(\koopvar)\\  
     1-\Sbar(1-\koopvar)
 \end{pmatrix}.
\end{equation}
So $S$ determined by \pccap is indeed is a linear map depending non-linearly on the input $k$, just as in the case of the orthonormalization procedure (\eqref{eqn:linearnonlinear}), and the action on the first component is equivalently given by the affine-linear map $\Sbar$. 

Using $\Sbar$ to learn only the first component of the $\chi$ function, we now can formulate the following explicit iterative \isokannI procedure
:


\begin{theorem}\label{thm:convergence}
    Let $\chi,\ \Sbar$ and $S$ be chosen as above in 
    \cref{eqn:chibar,eqn:shiftscale,eqn:fullS}.
    For generic $\bar\chi_0:\XX\rightarrow\RR$, i.e. containing components of $\ef_1$ and $\ef_2$, the \isokannI iteration  
    \begin{equation}\label{eqn:isokann1d}
        \bar\chi_{n+1} = \Sbar(\KK \bar\chi_n)
    \end{equation}
    converges to the $\chi$-function 
    \begin{equation}
        \label{eqn:chibarlim} \bar\chi = \lim_{n\rightarrow\infty} \bar\chi_n = \alpha \ef_1 + \beta \ef_2,
    \end{equation}
    for some $\alpha,\beta \in \RR$, which in turn solves the \isokann problem
    \begin{equation}\label{eqn:fixpointthm}\bar\chi = S\KK\bar\chi.\end{equation} 
\end{theorem}
\begin{proof}
Noting that $\Sbar$ is merely a shift-scale, i.e. affine linear, and $\KK\one = \one$ one obtains
\begin{align} \label{eqn:equalizer}
    (\Sbar\circ\KK)^n = \Sbar \circ \KK^n.
\end{align}
Looking at the eigendecomposition of $\bar\chi_0 = \sum_{i=1}^\infty a_i v_i$, we have
\begin{align}
    \KK^n \bar\chi_0 = \sum_{i=1}^\infty a_i \lambda_i^n \ef_i = a_1 \one + \sum_{i=2}^\infty a_i \lambda_i^n \ef_i
\end{align}
For large $n$, the contribution of the faster eigenfunctions $\ef_i,\,i > 2$ decays exponentially faster than the contribution of $\ef_1, \ef_2$. Hence the shift-scale is dominated by these slow eigenfunctions and with \eqref{eqn:equalizer} we have for some $\alpha,\beta \in\RR$
\begin{align}
    \bar\chi 
    = \lim_{n\rightarrow\infty} \left(\Sbar \circ \KK\right)^n \bar\chi_0
    = \lim_{n\rightarrow\infty} \Sbar \circ \KK^n  \bar\chi_0
    = \alpha \ef_1 + \beta \ef_2,
\end{align}
which proves \eqref{eqn:chibarlim}.
Noting that $\bar\chi(\XX)=[0,1]$ and $\KK$ acts as a shift-scale, and $\Sbar$ subsequently rescales $\KK\bar\chi$ back to the interval $[0,1]$ and we have $\Sbar(\KK \bar\chi) = \bar\chi$ which, given the implicit construction of $S$ from $\Sbar$ in \cref{eqn:fullS} shows the fixed-point result in \eqref{eqn:fixpointthm}.
\end{proof}

To summarize, we have shown how the representation of \isokann for $d=2$ as a scalar problem reproduces the classical \isokannI procedure \cite{Rabben2020isokann} where the role of the linear map $S$ (determined by \pccap) gets replaced by the (explicitly given) affine-linear $\Sbar$. This simpler representation allowed us to show convergence to a $\chi$-function and thus solving the $\isokann$ problem.

\subsection{Restoring the eigenfunctions}\label{sec:eigenfunction}
At the beginning of the article we intended to  compute the eigenfunctions of $\KK$. Whereas the $\chi$ functions only span the corresponding invariant subspace, we now show how we can restore the eigenfunctions from the solution $\chi$  to the \isokann problem $S \KK\chi = \chi$. This is equivalent to 
\begin{equation}
    \KK\chi = S^{-1}\chi
\end{equation}
which means that the action of $\KK$ on the subspace $\chi$ is given by the matrix $S^{-1}$. By means of a basis transformation (making use of the \emph{Moore–Penrose pseudoinverse} $\chi^+$) we can recover the eigenfunctions of $\KK$ from $S$ and $\chi$.

\begin{proposition}\label{prop:eigenfunction}
    Let $\chi=(\chi_i)_{i=1}^d$ be a column vector of $\chi_i$ component functions and $S$ be a full rank matrix such that 
    \begin{equation}
        S\KK\chi = \chi
    \end{equation}
    If $(X,\,\Lambda)$ is an eigendecomposition of $Q := \chi^+S^{-1}\chi$, i.e., 
    \begin{equation} \label{eqn:chieigenprob}
        QX = X\Lambda,
    \end{equation}
    then $E := \chi X$ are the eigenvectors of $\KK$ with eigenvalues $\Lambda$.
\end{proposition}
\begin{proof}
    By assumption we have $\KK \chi = S^{-1} \chi$. Inserting this into \eqref{eqn:chieigenprob}, multiplying with $\chi$ from the left and noting that $\chi^+\chi=\Id$ by definition, we arrive at
    \begin{equation}
         \KK \chi X = \chi X\Lambda
    \end{equation}
    which gives the desired result.
\end{proof}

\subsection{Computational procedure}
\label{sec:computation}

With these theoretical considerations, let us now discuss how to apply the algorithm in practice, i.e. using neural networks as function approximators and Monte Carlo  (MC) simulations for the Koopman evaluations.

To this end let us recall the main formula for the iterative update \eqref{eqn:isokann}:
\begin{equation}\label{eqn:isokanncopy}
    \chi_{n+1}(x_m) \leftarrow S_n \KKt \chi_n(x_m) = S_n \EE \left [\chi_n(X_\lag) \mid X_0 = x_m\right]
\end{equation}
Here we replaced the equality by a $\leftarrow$ which indicates classical supervised learning (with the common mean squared error loss) along multiple (possibly random) training points $x_m$, $m=1,...,\Nx$.
For a procedural description see \Cref{alg:isokann}.
The main challenges posed by this iterative scheme consist of (a) a representation of the function(s) $\chi_n$ and (b) the evaluation of the right hand side, i.e. the computation of $S_n$ and the evaluation of the Koopman operator.

For (a), the representation of the $\chi_n$, we chose neural networks as they promise good approximation properties in high dimensions and their differentiability will prove crucial for the following optimal control part.
In general any feed-forward architecture should be suitable and whilst convolutional networks could be especially suited due to the spatial structure of the state space $\XX$, in the example we will confine ourselves to a fully connected architecture for simplicity.
In any case, the update step for $\chi_{n+1}$ consists of a classical supervised learning routine with the labeled data
\begin{equation}
    D_n = \left\{ (x_m, s_m) \right\}, \quad s_m = S_n\KK \chi_n (x_m)
\end{equation}
generated by evaluation of the current $\chi_n$.
Because the $\chi_n$ are not expected to be changing a lot between the iterations it makes sense to initialize $\chi_{n+1}$ with the weights from $\chi_n$ as to transfer the already learned structure and speed up the learning. Note here that whilst we talk about different networks $\chi_n$ for each iteration $n$ to emphasize the iterative nature of \eqref{eqn:isokanncopy}, in practice we can update a single instance of the network. 

In this view, the learning procedure can be seen as iterative supervised batch learning, where the whole data batch $D_n$ is generated along a set of points $\{x_m\}$ using the current representation $\chi_n$.
The update step itself can be performed using any stochastic optimizer such as classical stochastic gradient descent or ADAM to minimize the empirical $L^2$ error 
\begin{equation}
\label{eqn:trainingloss}
    \min_{\chi_{n+1}} ! \sum_{(x_m, s_m)\in D_n} (\chi_{n+1}(x_m) - s_m)^2.
\end{equation}

Considering (b), the evaluation of the right hand side, let us start with the approximation of the Koopman operator. For a given training point $x_m$, we use its representation as an expectation value and approximate the action of the Koopman operator by a Monte-Carlo sum. Each simulation consists of starting $\Nmc>0$ trajectories at the point $x_m$ and propagating them according to the SDE \eqref{eq: not controlled sde} using an SDE integrator, such as the Euler-Maruyama scheme, for the lag time $\lag$ and storing their end points $y_{k,m}, k=1,...,\Nmc$. The action of the Koopman operator at $x_m$ is then approximated by the empirical average
\begin{equation}
\label{eq:koopman-chi}
    \koopvar_m := \KK \chi_n (x_m) = \EE \left [\chi_n(X_\lag) \mid X_0 = x_m\right] \approx \frac{1}{\Nmc} \sum_{k=1}^{\Nmc} \chi_n (y_{k,m}).
\end{equation}
To summarize, for each $x_m$ we average the evaluation of $\chi_n$ at $\Nmc$ propagated positions obtained by SDE simulations.

What now remains is the application of the shift-scale $S_n$. In case of \isokannI the action of $S_n$ is determined by the shift-scale $\Sbar$ from \eqref{eqn:shiftscale} which depends on the (global) extrema of the input function $k$. In practice we thus use the empirical extrema over the observed data $\{\koopvar_m\}$ to directly compute $s_m = \left(\Sbar(\koopvar)\right)_m$ without explicitly constructing $S_n$.
In higher dimensions we compute the matrix $S_n$ using \pccap to find the transformation that maps the columns of the matrix $K = [\koopvar_1\, ...\, \koopvar_m]$ into the unit simplex.
Note that the use of the empirical extrema requires that the training points $x_m$ indeed cover the areas where $\KK \chi_n$ becomes (approximately) extremal.

This brings us to the choice of the $\Nx$ training points $x_m$. In principle \isokann can be applied to find the $\chi$ functions of a system based on a fixed set of precomputed or assimilated trajectories (replacing the SDE integration). However its iterative nature makes it especially useful in the synthetic data regime where the trajectories are computed on-line as it allows adapting the training points $x_m$, and as we will see the trajectory simulations too, to the information obtained so far.

Since the $\chi$-functions can be interpreted as reaction coordinates \cite{Ernst2019} we suggest \emph{"$\chi$-stratified" sampling} of the $x_m$, i.e. such that $\chi(x_m)$ is approximately uniform in $[0,1]$.
In practice we achieve this by subsampling from the pool of start and end points of the previous simulations, 
\begin{equation}
    P_n=\left(\bigcup_m x_m^{(n-1)} \right)\cup\left(\bigcup_{m,k}y_{k,m}^{(n-1)}\right).
\end{equation}
We then draw $M$ stratified uniform samples $u_i\in[0,1]$ by sampling uniformly from each of $M$ equally sized partitions of the interval $[0,1]$.
Finally we chose those $p_i \in P_n$ such that  $\chi_n(p_i)$ is the closest to one of the $u_i$.
We furthermore retain those samples which were extremal in $\chi_n$ to facilitate good approximation of the extrema by $\Sbar$ or \pccap.

Heuristically speaking we obtain samples $x_m$ which are uniform in $\chi$ and hence provide good coverage or "bridges" along the transition region, thus facilitating an efficient "flow of information" during the power-iteration process.
Furthermore the regions with a higher variation of $\chi$, i.e. those that are "harder to learn", will also obtain more samples which in turn is also beneficial for the training of the neural network itself.


Last but not least the samples obtained this way followed the system's ergodic dynamics and will therefore approximate the stationary distribution (restricted on each level set of $\chi$).
This allows us to evade the curse of dimensionality by restricting the sampling to physically meaningful samples along the reaction paths.

To summarize, \emph{$\chi$-stratified sampling} allows us to sample uniform along $\chi$ and stationary conditioned on it without much additional cost and adapted to the learning process.

The main \isokann routine can be summarized by three loops ($\Niter$ power iterations, $\Nx$ training points, $\Nmc$ trajectories) which can be loosely tied to the three main ingredients of \isokann:

\begin{enumerate}
    \item[(1)] The power iteration learning the dominant subspace.
    \item[(2)] The neural network approximation of the next $\chi$ iterate.
    \item[(3)] The Monte-Carlo simulation of the Koopman evaluation.
\end{enumerate}

In the outermost loop (1) we perform the power iteration $\chi_n$ to $\chi_{n+1}$. Whereas in theory we have convergence for $n\rightarrow\infty$, we chose to terminate after a fixed number of iterations $\Niter$. This could be replaced by classical convergence criteria such as relative and absolute tolerances. Note that the rate of convergence and hence the required number of power iterates $\Niter$ depends on the eigenvalues of the Koopman operator where a bigger spectral gap implies faster decay of the non-dominant spectrum and hence faster convergence.

When training the neural network (2) we loop over a batch of labeled training data at the $\Nx$ training points $x_m$ which in turn (3) require $\Nmc$ individual trajectory simulations.
Both the number of training points $\Nx$ as well as trajectories per point $\Nmc$ depend on the step sizes chosen for the neural network optimizer. Since in practice the evaluation of the Koopman operator is rather expensive compared to the neural network update, it may be efficient to perform multiple update steps on the same batch of data before proceeding to the next iteration.

Note that the variance of the training data, scaling with $\Nmc^{-1}$, is particularly high for metastable systems due to the impact of rare transitions. Whereas above we proposed the use  \emph{$\chi$-stratified} subsampling as a heuristic to deal with sampling in $\XX$ space, we will now address the problem of variance in the "$\KK\chi_n$-direction" using the techniques of optimal control and importance sampling.

\begin{algorithm}[H] 
\caption{\isokann}\label{alg:isokann}
\begin{algorithmic}[1]
\Require $\Niter$ the number of power iterations\\
$\Nx$ the number of $x$-samples\\
$\Nmc$ the number of Koopman Monte-Carlo samples\\
$\chi_0$ initial neural network 
\Ensure $\chi_{\Niter}$ approximates a $\chi$ function (c.f. \cref{eq:isofix})\\
\For {$n=0$ to $\Niter-1$} 
    \For{$m=1$ to $\Nx$} 
        \State $x_m \gets$ \Call{sample$X_0$}{}$() $ \Comment{sample training points}
        \For {$k=1$ \textbf{to} $\Nmc$} 
            \State $y_{k} \gets$ \Call{sample$X_T$}{$x_m$} \label{algline:sde} \Comment{simulate trajectories}
        \EndFor
        \State $\koopvar_m \gets \frac{1}{\Nmc} \sum_k \chi_{n}(y_k) $ \label{algline:koopman}\Comment{Koopman approximation}
    \EndFor
    \State $s \gets \Call{S}{\koopvar}$  \Comment{transformed target data}
    \State $\Delta \gets \nabla_{\theta_n} \sum_m (\chi_{n}(x_m) - s_m) ^2$
    \Comment{compute loss gradient}
    \State $\chi_{n+1} \gets \Call{optim}{\chi_n, \Delta}$ \Comment{train the neural network}
\EndFor
\State \Return $\chi_{\Niter}$
\Statex
\Statex \hspace{-1.7em}\textbf{Subroutines}:\\
\Call{sample$X_0$}{}: subroutine sampling the starting points $x_m$ 
- either uniform or \emph{$\chi$-stratified} (\Cref{sec:computation}). \\
\Call{sample$X_T$}{}: SDE solver, e.g. Euler-Maruyama
- either uncontrolled or controlled (\Cref{sec:controlapplication}).\\
\Call{S}{}: empirical shift-scale or \pccap (\Cref{sec:isokannI} or \ref{sec:transformationS})\\
\Call{optim}{}: gradient based optimization of the neural network (e.g. 100 SGD steps)
\end{algorithmic}
\end{algorithm}

\section{Optimal sampling of Koopman eigenfunctions and \texorpdfstring{$\chi$}{chi}-functions}\label{sec:optcontrol}
In this chapter we first recall importance sampling, before showing how the theory allows to better sample eigenfunctions of the Koopman operator and \isokann $\chi$-functions. We conclude this chapter with a numerical example.

\subsection{Importance Sampling for random variables}\label{sec:oisRV}
Importance sampling allows to express the expectation value of an observable $f>0$ with respect to some distribution $p$ by an expectation value with respect to some other distribution $q$ (with $p \ll q$) by the formula
\begin{align}
    Z := \EE_p[f] = \EE_q\left[f \frac{\dd p}{\dd q}\right],
\end{align}
where the observable $f$ is reweighted by the Radon-Nikodym derivative $\frac{\dd p}{\dd q}$. 
It is easy to see that by choosing $q^*$ such that $\frac{\dd p}{\dd q^*} = \frac{Z}{f}$ (i.e. $q^* = \frac{f}{Z} p$), we have 
\begin{align}
    Z = \EE_{q^*}\left[f \frac{Z}{f}\right] =\EE_{q^*}[Z].
\end{align}
Now, since $Z$ is a constant, it can be computed with a single (reweighted) $f$-sample from $q^*$. 
We therefore refer to this importance sampler with sampling distribution $q^*$ as \emph{zero-variance-sampler}, or \emph{optimal-importance-sampler}.
Note however that we needed to know the (a priori unknown) result $Z$ in order to define the optimal sampling distribution $q^*$.

\subsection{Optimal Importance Sampling for Diffusion Processes}\label{sec:oisdiff}
Since we are working with diffusion processes, importance sampling is further complicated in that the measures are path measures and admit no probability density function. However, importance sampling can still be generalized to stochastic processes. Let us first consider the diffusion process of Eq.~\eqref{eq: not controlled sde}.
In general, one is interested in computing the expectation of path-dependent quantities. Let us define the \emph{work} along a trajectory over $[0,T]$ by the accumulation of a \emph{running cost} $f$ and a \emph{terminal cost} $g$ as:
\begin{equation}
    \label{eq: work}
    W_{t, T}(X) \coloneqq \int_t^T f(X_s, s) \dd s + g(X_T) .
\end{equation}
One then is interested in estimating expectation values of the form
\begin{equation}
    \label{eq: psi}
    \psi(x, t) \coloneqq \EE_{X_t=x} \left[ \exp\left(- W_{t, T}(X)\right) \right] = \EE_{X_t=x} \left[ \exp \left(-\int_t^T f(X_s, s) \dd s - g(X_T) \right)\right] .
\end{equation}

Girsanov's theorem builds the bridge from importance sampling to diffusion processes by allowing to sample from another diffusion processes. In particularly, it allows us to compute the change of measure in terms of the Radon-Nikodym derivative. 
To this end let us introduce the controlled process ${X^u = (X_t^u)_{t \geq 0}}$
\begin{equation}
\label{eq: controlled sde}
\dd X_t^u = (b(X_t^u) + \sigma u(X_t^u, t)) \dd t + \sigma \dd B_t,
\end{equation}
with an admissible control term $u$ acting as an external forcing to the original dynamics. Note that with zero control $u=0$ one recovers the original dynamics $X = X^{u=0}$. 
Let $\mathcal{P}$ denote the path measure induced by $X$ and $\mathcal{Q}$ denote the measure induced by $X^u$. 
According to Girsanov's theorem the change of measure from $\mathcal{Q}$ to $\mathcal{P}$ (analogous to $\frac{\dd p}{\dd q}$ above) along a given controlled trajectory $X^u$ is then given by
\begin{equation}\label{eqn:girsanov}
    G_{t, T}(X^u) \coloneqq\frac{\dd\mathcal{P}}{\dd\mathcal{Q}} \big |_{[t, T]} (X^u) = \exp \left(- \int_t^T u(X_s^u, s) \cdot \dd B_s  - \frac{1}{2} \int_t^T \left| u(X_s^u, s) \right| ^2 \dd s \right)
\end{equation}
which in turn provides an unbiased estimator of $\psi$ in terms of the controlled process:
\begin{equation}
\label{eq: girsanov unbiased estimator}
\psi(x,0) = \EE_{X_0 = x} \left[ \exp{\left( -W_{0, T}(X) \right)}\right] = \EE_{X_0 = x} \left[ \exp{\left( -W_{0, T}(X^u) \right)} G_{0, T}(X^u)\right] .
\end{equation}

Note that even though the expected value of this estimator is the same for any control $u$ its variance will vary.
Analogous to the case above there exists an optimal measure corresponding to an optimal control $u^*$ for which the controlled estimator exhibits zero variance \cite{Hartmann2017,Nusken2021solving}:
  
\begin{theorem}\label{thm:optcontrol}
    The optimal control $u^*$ is given by
    \begin{align}
        u^*(x,t) 
        = \sigma^\top \nabla_x \log \psi(x,t)
    \end{align}
    and leads to the zero variance estimator for $\psi$
    \begin{align}
        \psi(x,0) 
        = \EE_{X_0=x}\left[ \exp(-W_{0,T}(X)) \right] 
        \overset{a.s.}{=} G_{0,T}(X^{u^*}) \exp(-W_{0,T}(X^{u^*})) \text{ with }X^{u^*}_0=x.
    \end{align}
\end{theorem}

\subsection{Optimal sampling of eigenfunctions of the Koopman operator}\label{sec:oiskoop}

We will now show how this optimal control theorem can be used to evaluate the Koopman operator.

\begin{corollary}\label{cor:optcontrol}
Let $h \in L^\infty(\XX)$ be a function. A single realization of the controlled process $X^u$ starting in $X^u_0 = x$ with  control
\begin{equation} \label{eqn:optcontrolK}
    u(x,t) = \sigma^\top \nabla_x \log (\KK^{T-t} h)(x)
\end{equation}
then gives the evaluation of $\KK^T h$ at that point $x$:
\begin{equation}
    (\KK^T h )(x) \overset{a.s.}{=} G_{0,T}(X^u) h(X^u_T).
\end{equation}
\end{corollary}
\begin{proof}
With $f(x) = 0$ and $g(x)=-\log h(x)$ \cref{eq: psi} becomes
\begin{equation}
    \psi(x,t) 
    = \EE_{X_t=x} \left[h(X_T)\right] 
    = \EE_{X_0=x} \left[h(X_{T-t})\right] 
    = (\KK^{T-t}h)(x),
\end{equation}
Application of \Cref{thm:optcontrol} then leads to the desired result.
\end{proof}

This result shows us that in order to compute the optimal control for evaluating $\KK h$ we need to have access to the derivatives of $\KK h$. This conundrum is in line with the general optimal importance result and comes at no surprise. In the next section we will argue how this result can still be of use for \isokann where the convergence of the $\chi_n$ provides us with an approximate description which we will use to compute the control.

Let us for now consider the case where the observable of interest $h$ is an eigenfunction of $\KKt$, i.e. $h = \ef_i$ with eigenvalue $\lambda_i(\lag)$. In this case we can replace the action of the Koopman operator by its eigenvalue, which in turns cancels out after application of $\nabla \log$ and results in a time-independent control:

\begin{align}
    u^*(x,t) 
    = \sigma^\top \nabla_x \log (K^{T-t} \ef_i)(x) 
    = \sigma^\top \nabla_x \log (\lambda_i(T-t) \ef_i(x))  
    = \sigma^\top \frac{\nabla_x \ef_i(x)}{\ef_i(x)}.
\end{align}

This simple example provides a good point to get a feeling for how optimal importance sampling works. We can see that the control pushes the system in the direction of (relative) maximal ascent. If the system follows the forcing its expected evaluation increases, but the path taken also has an increased probability, resulting in a decreasing Girsanov weight.
Both increments happen on a commensurate "relative scale": the expectation value gets pushed by an amount relative to its current value ($\frac{\nabla v}{v}$) whereas the reweighting is adjusted relatively in magnitude due to the exponentiated integral (which may be seen as an infinite product over all time-points) in \eqref{eqn:girsanov}.
In this way the increase of the observable is balanced with the decreasing weight exactly so that no matter the path taken these always equalize and one obtains a zero-variance sampler.

Unfortunately however, the control becomes singular whenever $\ef_i(x) = 0$. According to the Perron-Frobenius theorem, every non-trivial eigenfunction is unsigned, i.e. crosses the 0 at some point, so we have to find a way around that problem.
We can alleviate this problem by shifting and (anticipating the form of $\Sbar$ in \eqref{eqn:shiftscale}) also rescaling the eigenfunction.

Denote the shift-scaled eigenfunction by 
\begin{align}
    \chi &:= \alpha \ef_i + \beta \one, \quad \alpha, \beta \in \RR, \quad \chi>0 .
\end{align}
Due to linearity of $\KKt$, we obtain
\begin{align}\label{eqn:chilincomb}
    (\KK^{\lag} \chi)(x) =  (\KK^{\lag} \alpha \ef_i)(x) + \beta = \alpha \lambda_i(\lag) \ef_i(x) + \beta
\end{align}
and thus after application of \Cref{cor:optcontrol} the control in terms of $\chi$ is

\begin{equation}
\label{eqn:optchicontrol}
    u^*(x,t) 
    = \sigma^\top \nabla_x \log (\KK^{T - t} \chi)(x)
    = \sigma^\top \frac{\nabla \chi(x)}{\chi(x) + \frac{\beta}{\lambda_i(T-t)} - \beta} .
\end{equation}

In this case we see that control is time dependent (which makes sense as the relative contributions of the dominant eigenfunctions to the expectation value change over time). From $\lambda_i(t) \le \lambda_i(0) = 1$ we can conclude 
that $u^*(\cdot, t) = \gamma \sigma^\top \frac{\nabla\chi}{\chi}$ with function $\gamma$ monotically increasing with $\gamma(1) = 1$, i.e. the control is pointing in the same direction as for the pure eigenfunction case, starting weaker and increasing until hitting the full magnitude at $t=T$.

Note that the requirement for $\chi$ functions to satisfy $0 \le \chi \le 1$, which so far was merely motivated by their interpretation as macrostates, now also facilitates their optimally controlled importance sampling.

\subsection{Application to \isokann}\label{sec:oisiso}
\label{sec:controlapplication}
In the previous section we have shown how to obtain a zero-variance sampler for the Koopman operator in terms of the gradient of its solution, either for general observables $h$ or (shift-scaled) eigenfunctions $\ef_i$ ($\chi$). In either case the solution has to be known a priori as to compute the control. We will now argue how to integrate this result into the \isokann procedure.

The main idea of using optimal importance sampling in \isokann is to use the intermediate results $\chi_n$ and $S_n$ to compute a \emph{pseudo-optimal control} as to lower the sampling variance.
We therefore have to assume that using an approximation to the optimal control indeed leads to a variance reduction. Whereas we do not know of any proof to this statement it was shown  that the objective of the associated optimal control problem is indeed convex in the control \cite{Lie2016convex} which leads us to conjecture that that the variance should be well-behaving for approximate optimal controls as well.

Note that the importance sampler \eqref{eq: girsanov unbiased estimator} is unbiased for any control. Thus even if the above assumption does not hold \isokann would still converge, albeit slower, as long as the variance does not become unbounded\footnote{Which could always be ensured by e.g. clipping the control, thus bounding the Girsanov reweighting term and in turn also the overall sampling variance}, under the usual conditions for stochastic gradient descent convergence (i.e. decaying learn rate).

Let us recall the equation for the control \eqref{eqn:optcontrolK}:
\begin{equation}
    u(x,t) = \sigma^\top \nabla_x \log (\KK^{T-t} h)(x)
\end{equation}
In order to compute the differential of $\KK^{T-t}$ at $\chi_n$ we assume sufficient convergence of \isokann together with \eqref{eqn:isokann},
\begin{equation}
\chi_n \approx \chi_{n-1}, \quad \chi_n = S_{n-1}\KKt \chi_{n-1},
\end{equation}
to approximate the action of $\KK^\lag$ by $(S_{n-1})^{-1}$:
\begin{equation}
    (S_{n-1})^{-1} \chi_n = \KKt \chi_{n-1} \approx \KKt \chi_n.
\end{equation}
Using the semi-group property of $\KK$ and the matrix logarithm we can extend this to other lag times $T-t$ to obtain the matrix approximation $\widetilde \KK^{T-t}$
\begin{equation} \label{eqn:Ktilde}
\KK^{T-t}\chi_n \approx \widetilde \KK^{T-t} \chi_n := \exp\left( \frac{T-t}{\lag} \log \left((S_{n-1})^{-1}\right)\right) \chi_n.
\end{equation}

Note that in the general $d$-dimensional case the expectation values in the \isokann iterations $\KKt \chi_n$ are vector valued. Optimal importance sampling however works only in the scalar case. Therefore we have to compute an individual control for sampling each component $(\KKt \chi_n)_i$ individually\footnote{In conjunction with \emph{$\chi$-stratified} sampling this results in multiple search directions, each exploiting the assumed location of one of the metastabilities. Since this also implies moving away from the respective other metastabilties (and beyond the current one) we have hopes that this interplay between the search directions may automatically provide a balance between exploration and exploitation.}.

Thus using the optimal control \Cref{cor:optcontrol} together with the matrix approximation \eqref{eqn:Ktilde} we can compute the pseudo-optimal control for the $i$-th component of $\KKt \chi_n$ explicitly
\begin{align}
    u^*_i(x,t) = \sigma^\top \nabla_x \log \left(\sum_j \widetilde \KK^{T-t}_{ij}\chi_j(x)\right)=\sigma^\top \frac{\sum_j \widetilde \KK^{T-t}_{ij} \nabla_x \chi_j(x)}{\sum_j \widetilde \KK^{T-t}_{ij}\chi_j(x)}
\end{align}

In the case of \isokannI the action of the Koopman operator on $\chi_n$ converges to a shift-scale as in \eqref{eqn:chilincomb} and we can therefore estimate the parameters $\alpha$, $\beta$ and $\lambda_2$ from the extrema of $\KK \chi_{n-1}$ as to apply the explicit control for the shift-scaled eigenfunction \eqref{eqn:optchicontrol}.

Now that we know how to compute the control we can modify the algorithm (\Cref{algline:sde}) by sampling the trajectories according to the controlled SDE \eqref{eq: controlled sde}. In order to compute the reweighting \eqref{eqn:girsanov} we have to either save the trajectory and noise, or integrating it on the fly in an addition SDE component with
\begin{equation}
    G = \exp(-g_T)\quad 
    \dd g_t = \frac{1}{2}u(X^u_t,t)^2 \dd t + u(X^u_t,t) \cdot \dd B_t \quad 
    \dd g_0 = 0.
\end{equation}
Finally, for the Koopman Monte Carlo approximation (\Cref{algline:koopman}) we average over the $\chi$ evaluations at the endpoints of $\Nmc$ independent trajectories starting in $x_m$ weighted  with their respective weights $G_m$:
\begin{equation}
\KK \chi_n (x_m) \approx \frac{1}{\Nmc}\sum_{k} \chi_n(X^u_{T,m}) G_m.
\end{equation}

In this way (and with the above assumption) we obtain a feedback loop where better approximation of the $\chi$ functions results in in a better approximation of the action of $\KK$ and hence in a better approximation of the optimal control. This pseudo-optimal control in turn decreases the sampling variance which facilitates better approximation of the power iterates, i.e. the $\chi$ function.

As a proof of concept we will now illustrate the reduction of variance at the hand of the classic double-well potential.

\subsection{Example: Controlled \isokannI for the double well}\label{sec:oisexample}

Let us consider the controlled process as stated in \eqref{eq: controlled sde} for the double-well potential \eqref{eqn:doublewell}
which leads to the simplest problem exhibiting metastable behavior and hence a challenging sample variance. 
In our experiments we compare the training performance of the \isokann algorithm, both, with and without the control \eqref{eqn:optchicontrol}. 

We start with a randomly initialized fully connected network $\chi_0: \mathbb{R} \rightarrow \mathbb{R}$ with sigmoidal activation functions and 2 hidden layers, each of size 5 (i.e. with layer sizes $1\times5\times5\times1$).
For each network generation $n$, we compute Monte Carlo approximations of the Koopman expectation at $\Nx=30$ positions. These are initially drawn uniformly from the interval $[-2,2]$ and subsequently obtained by \emph{$\chi$-stratified} subsampling as described in \Cref{sec:computation}.
From each starting position we then simulate $\Nmc=20$ trajectories using the SROCK2 SDE integrator of strong order 1 with step-size $\Delta t=0.001$. The next generation $n+1$ is trained against these $\Nx$ training points by $\Nlearn=500$ stochastic gradient descent steps using the ADAM optimizer (with learning rate $\eta = 0.001$). We repeat this evaluation-training procedure (corresponding to a single power iteration) for a total of $\Niter=50$ iterations.

For each experiment we monitor the root of the training loss \eqref{eqn:trainingloss}, i.e. the root mean squared error, 
\begin{align}
\text{RMSE} &\coloneqq \sqrt{ \frac{1}{\Nx}\sum_m (\chi_{n}(x_m) - s_m) ^2} \\
\intertext{and the mean standard deviation of the MC estimator \eqref{eq:koopman-chi}}
\text{MSTD} &\coloneqq 
\frac{1}{\Nx} \sum_m
\sqrt{\frac{1}{\Nmc} \sum_k (\chi_{n}(y_{k,m}) - \mu_m ) ^2}, 
\quad \text{where} 
\quad \mu_m =\frac{1}{\Nmc} \sum_k (\chi_{n}(y_{k,m}))
\end{align}
over the training phase of $\Niter=50$ iterations with $\Nlearn=500$ training steps each.


\begin{figure}
\begin{subfigure}{.49\textwidth}
\includegraphics[width=\textwidth]{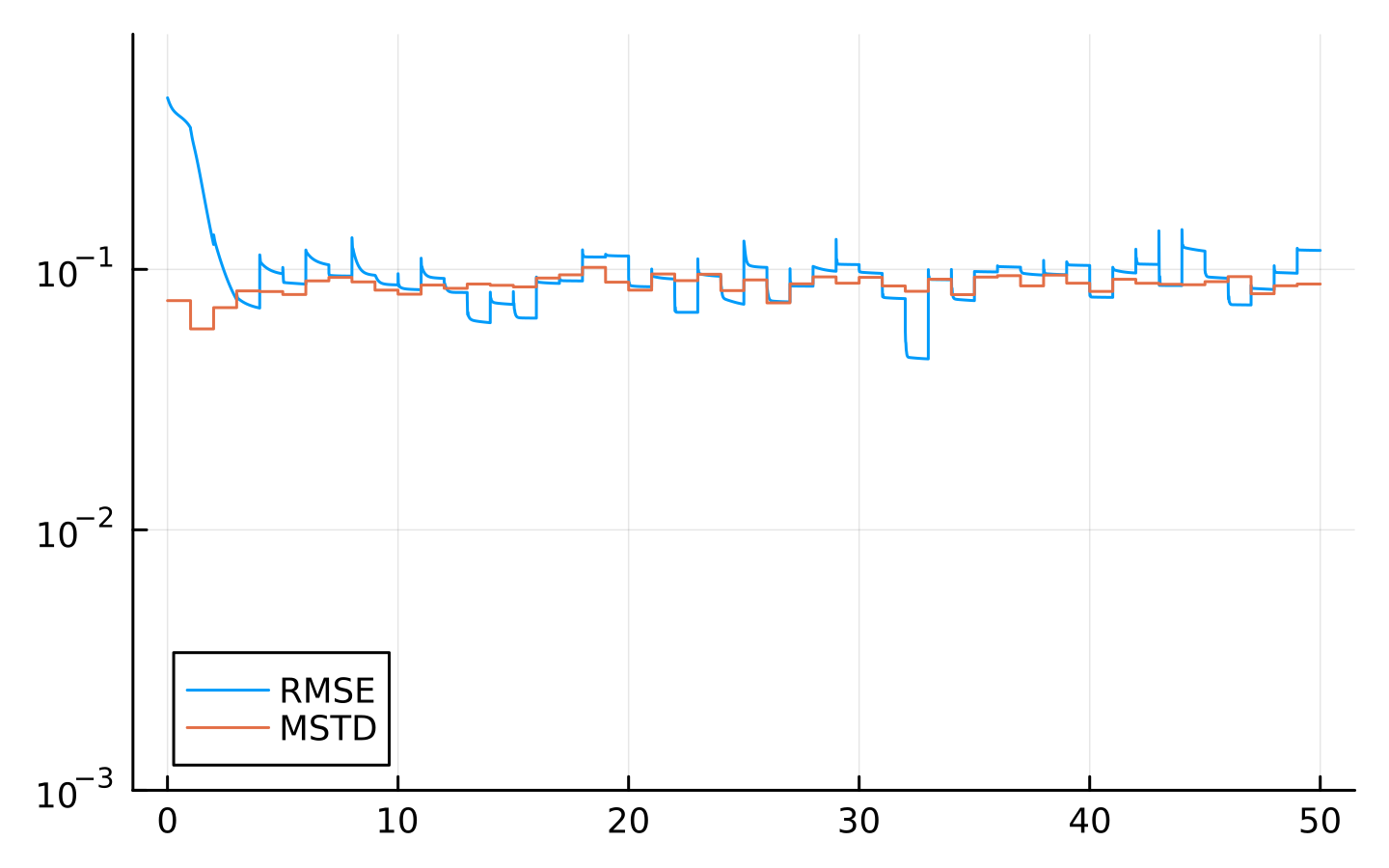}
\caption{without control}
\label{fig:loss-false}
\end{subfigure}
\begin{subfigure}{.49\textwidth}
\includegraphics[width=\textwidth]{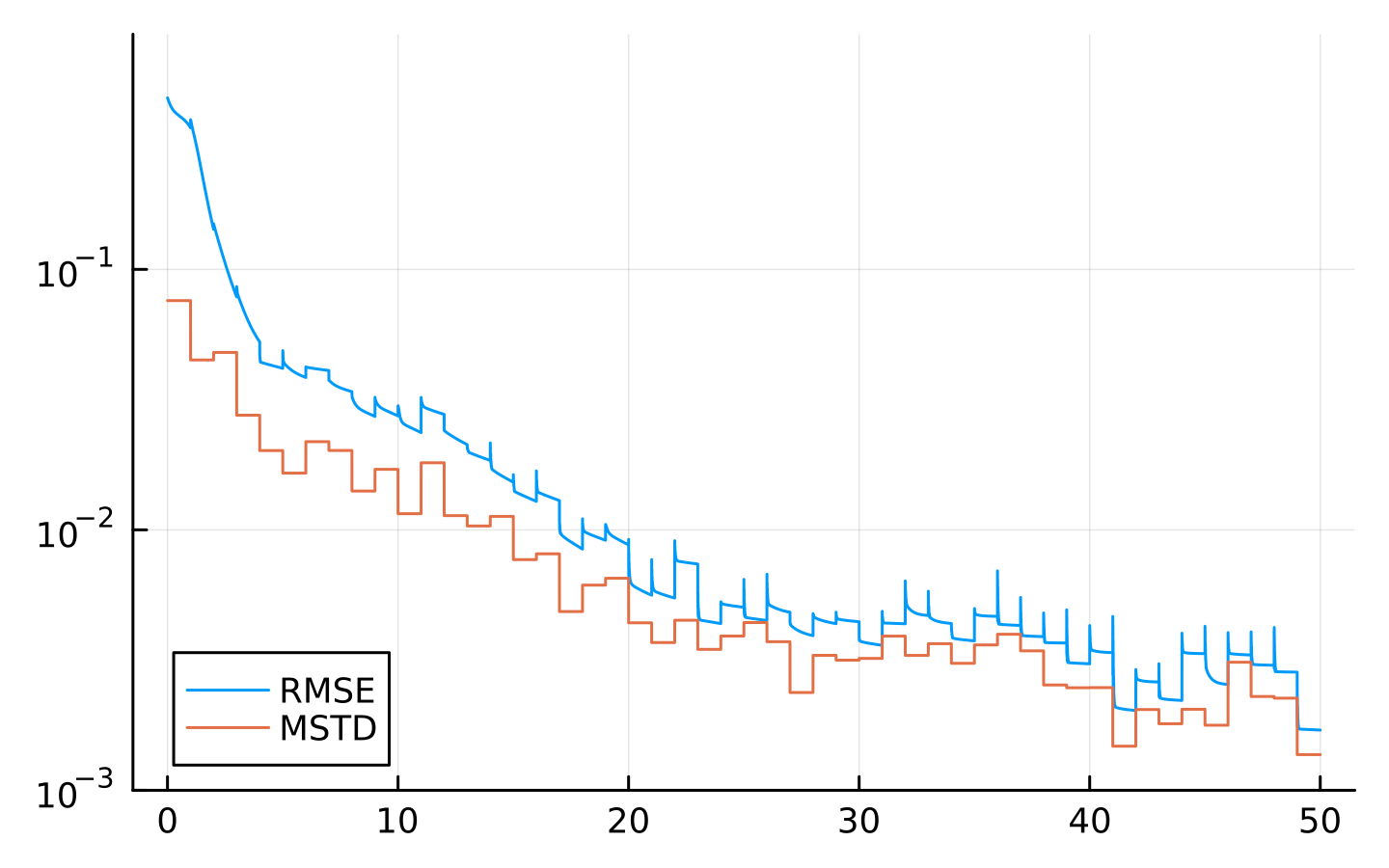}
\caption{with control}
\label{fig:loss-true}
\end{subfigure}
\caption{Training performance over $50$ power iterations / batches with 500 ADAM steps each. The blue line shows the training loss and the red line shows the standard deviation of the Monte-Carlo samples in the training data.
}
\label{fig:loss}
\end{figure}

Let us now compare the uncontrolled with the controlled experiment. \Cref{fig:loss} shows the (square root of) our training loss together with the standard deviation of the Monte Carlo estimator for the two cases of study. 
In \Cref{fig:loss-false} we observe that the uncontrolled system quickly (after 3 iterations) approaches its plateau at an error of about $10^{-1}$ but afterwards the training loss does not decrease any further. This comes at no surprise since the training data exhibits noise of the same magnitude, and we cannot expect the average loss to be lower then the noise in the data.
Note however, that even though the loss itself seems to have leveled off this does not necessarily mean that the solution does not improve: Whilst the empirical loss will necessarily remain at the level of the noise, the solution could still converge due to the inherent averaging of that noise in the SGD method.

Looking at the controlled experiment, we observe in \Cref{fig:loss-true} that for the loss behaves similar to the uncontrolled experiment for the first 3 iterations, reaching a value of $10^{-1}$. From there on however the loss decreases further getting close to $10^{-3}$ after $50$ power iterations and still not having hit a plateau. Notice that the training noise, in strong contrast to the uncontrolled case, decreases rapidly from the beginning of the training. It is furthermore interesting to see that the training loss seems to be following the noise level closely, indicating that the sampling variance is indeed of high importance.

Looking more closely at these performance plots we can furthermore identify the individual training batches: The MSTD is piecewise constant along each such batch, since the training data (and hence its standard deviation) is updated only inbetween the neural network training loops. Whereas the RMSE plateaus continiously during each batch, we can observe how it jumps up slightly at the beginning of each new batch. These jumps are caused by overfitting to the previous batch (the plateau) without generalization to the following batch (the flank)
\footnote{To increase the efficiency of the algorithm one could therefore take the plateau of each flank as indication to interrupt the current training and generate a new training batch.}.

\begin{figure}
\begin{subfigure}{.495\textwidth}
\includegraphics[width=\textwidth]{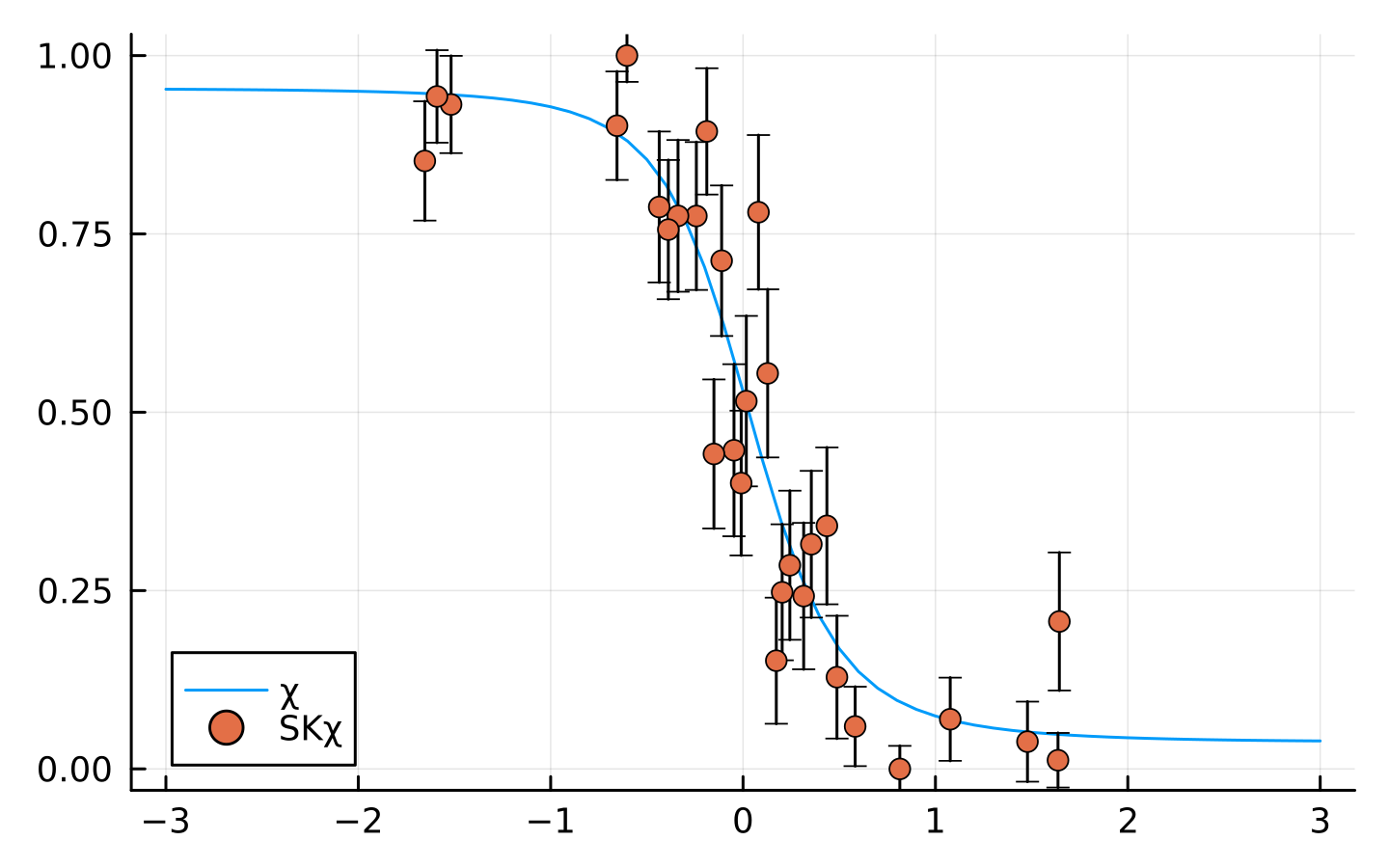}
\caption{without control}
\end{subfigure}
\begin{subfigure}{.495\textwidth}
\includegraphics[width=\textwidth]{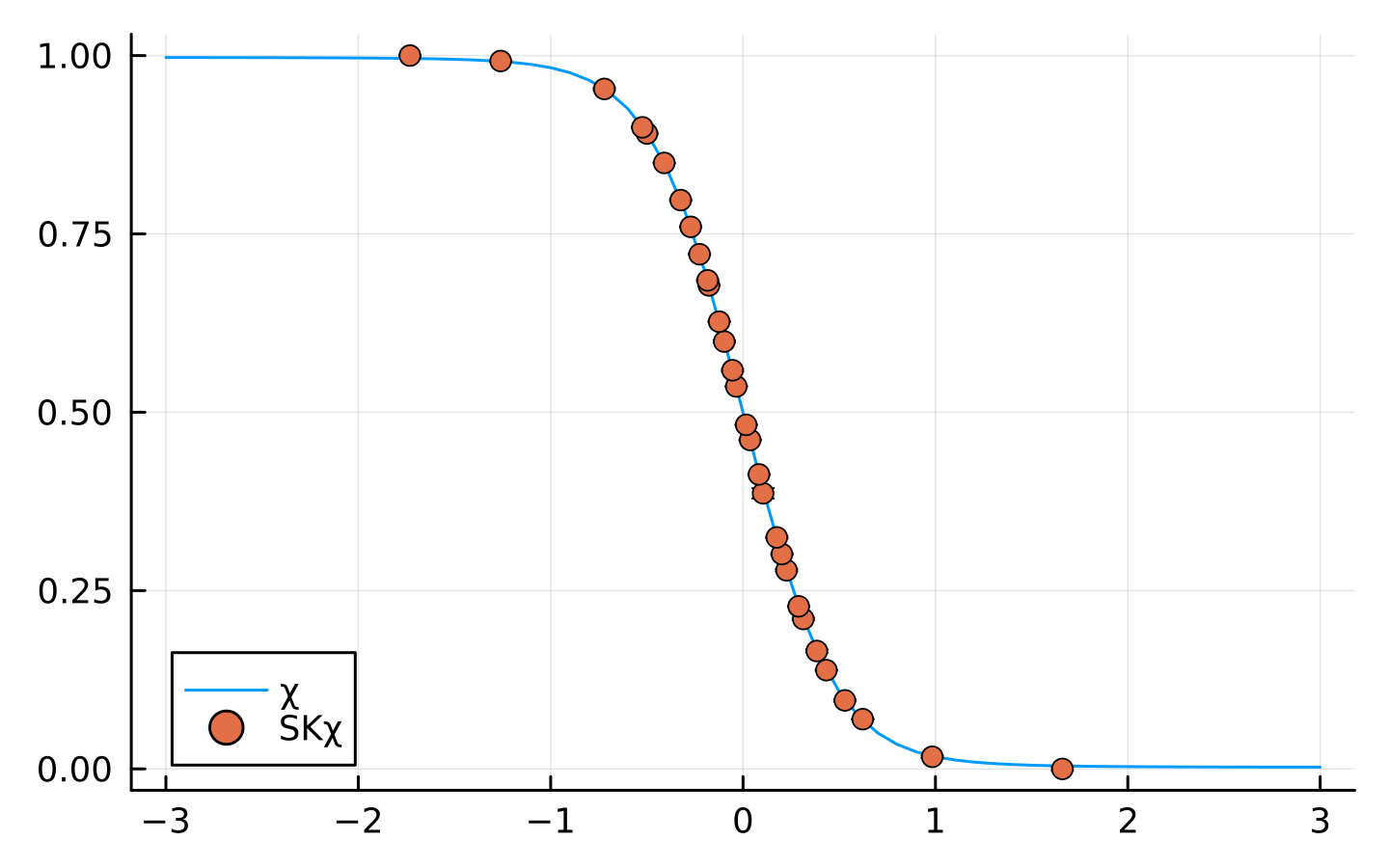}
\caption{with control}
\end{subfigure}
\caption{The blue line shows the learned $\chi$ function at the end of training. The red dots show the train target, i.e. the Koopman evaluations at the $x_m$, with the Monte-Carlo standard deviation as error bars.}
\label{fig:chi-fit}
\end{figure}

Let us finally look at \Cref{fig:chi-fit}, which shows the different learned $\chi$ functions after the application of the \isokann algorithm together with the evaluations of $SK \chi_{n-1}(x_m)$ at the $\Nx$ random locations $x_m$. The error bars represent the standard deviation of the individual Monte-Carlo estimators, i.e. the noise in the training data. We see that both learned $\chi$ functions qualitatively match the expectation. The uncontrolled case however has problems reaching $0$ resp. 1 at the boundaries, which can be understood as a result of the noise and the subsequent noisy estimation of the empirical shift-scale. Last we notice that the Monte-Carlo standard deviation for the Koopman evaluation at the $\chi$-sampled positions is considerably lower by the controlled approach.

\section{Conclusion}





In this article we started by enhancing \isokann by new theoretical results that prove the strengths of \isokann, namely a convergence proof (Thm. \ref{thm:convergence}) and a method for reconstructing eigenfunctions from \chifunctions (Prop. \ref{prop:eigenfunction}).  We also proposed a new adaptive sampling strategy, called $\chi$-stratified sampling, which complements \isokann well and deserves further investigation. 
Formulating \isokann in terms of the transformation $S$ \eqref{eqn:isokann}, we paved the way for higher-dimensional \chifunctions while generalizing the original \isokannI.
However, whereas we argued for using \pccap for the construction of $S$ for $d>1$ a more detailed study and a proof of convergence for this case remain open for future work.

The second main contribution in this article is the introduction of importance sampling into \isokann. Whereas we know that the resulting estimator is unbiased we argued only heuristically why the variance should not explode. A proof of convexity of the variance in the control, or even better, a proof of the convergence of control in \isokann is still missing.
Note furthermore, that the concept of optimal importance sampling may be useful for the iterative solution of Koopman evaluations in general (c.f. Cor. \ref{cor:optcontrol}).
 
An important next step would be to apply controlled \isokann to an actual molecular dynamics (MD) system as to test how well the introduced techniques fare with the complexities of real world problems.
This however requires a way to run many trajectories with different start locations and low-overhead as well as to inject the optimal control into the MD simulations. We hope that once these interfaces are implemented, \isokann will enhance the research of molecular systems.

\section*{Acknowledgement}
\theacknowledge

\section*{Code availability}
The code used for the numerical examples is available as a Julia package on GitHub at \\
\url{https://github.com/axsk/OptImpSampling.jl} with the tag \texttt{jmp}
\footnote{In order to reproduce the experiments and plots of this paper run \\
\texttt{
julia> Pkg.add("\url{https://github.com/axsk/OptImpSampling.jl}\#jmp"); \\
julia> import OptImpSampling; OptImpSampling.paperplots()}}.

\printbibliography
